\documentclass[a4paper,10pt,twoside]{article}
\usepackage{amsmath,latexsym,amsfonts,amssymb,amsthm,mathrsfs,longtable,lscape,dsfont,fancybox,yfonts}
\usepackage{amsthm}
\usepackage{amssymb}
\usepackage{amscd}
\usepackage{amsfonts}
\usepackage{amsbsy}
\usepackage{fancyhdr}
\usepackage{graphicx}
\usepackage[all]{xy}
\newtheorem {teo} {Theorem} [section]
\newtheorem{lem}[teo]{Lemma}
\newtheorem {prop} [teo]{Proposition}

\title{Regularization in the Restricted Four Body Problem}
\author{Jaime Burgos--Garc\'ia\\
Published in Aportaciones
Matemáticas\\ Memorias de la Sociedad Matemática Mexicana, 45, pp.1-13.
 \thanks{Departamento de Matem\'aticas
UAM--Iztapalapa. Av. San Rafael Atlixco 186, Col. Vicentina, C.P.
09340, M\'exico, D.F. e--mail: jbg84@xanum.uam.mx} \thanks{Research article}}
\date{}
\begin{document}

\maketitle

\begin{abstract}
 The restricted (equilateral) four-body problem consists of three bodies of masses $m_{1}$, $m_{2}$ and $m_{3}$ (called primaries) lying in a Lagrangian configuration of the three-body problem, i,e,. they remain fixed at the apices of an equilateral triangle in a rotating coordinate system. A massless fourth body moves under the
Newtonian gravitation law due to the three primaries; as in the restricted three-body problem the fourth mass does not affect the motion of the three primaries. In this paper we show a global regularization of binary collisions of the infinitesimal body with two of the primaries.\\ \\ \begin{center}
 {\bf \MakeUppercase Resumen\\}
\end{center}
El problema restringido de cuatro cuerpos equil\'atero consiste de tres masas puntuales $m_{1}$, $m_{2}$, $m_{3}$ (llamadas primarias) que permanecen a todo tiempo en una configuraci\'on Lagrangiana del problema de tres cuerpos, es decir; las masas permanecen fijas en los v\'ertices de un triangulo equil\'atero en un sistema rotatorio. Un cuarto cuerpo de masa infinitesimal se mueve bajo la ley de gravitaci\'on universal de Newton que ejercen las tres masas puntuales; como en el caso del problema restringido de tres cuerpos, la cuarta masa no afecta el movimiento de las tres primarias. El objetivo principal de este art\'iculo es mostrar una regularizaci\'on global de colisiones binarias de la masa infinitesimal con dos de las primarias. Al final se muestra una aplicaci\'on de este proceso de regularizaci\'on.

\end{abstract}


\newpage
\section{Introduction}
Few bodies problems have been studied for long time in celestial
mechanics,  either as simplified models of more complex planetary
systems or as benchmark models where new mathematical theories can
be tested. The three--body problem has been source of inspiration
and study in Celestial Mechanics since Newton and Euler. In recent
years it has been discovered multiple stellar systems such as double
stars, triples systems. The restricted three body problem (R3BP) has
demonstrated to be a good model of several systems in our solar
system such as the Sun-Jupiter-Asteroid system, and with less
accuracy the Sun-Earth-Moon system.  In analogy with the R3BP, in
this paper we study a restricted problem of four bodies consisting
of three primaries moving in circular orbits keeping an equilateral
triangle configuration and a massless particle moving under the
gravitational attraction of the primaries. In the following discussion we focus on the study of the regularizations of binary collisions of the infinitesimal body with two of the primaries by a simple method similar to Birkhoff's which permit us to study the dynamic of the equations when they present discontinuities . As an application of the transformed equations by the regularization process it can be shown that some families of periodic orbits end up in a homoclinic connection. This last
phenomenon can be dynamically explained by the so called ``blue sky
catastrophe'' termination, a rigorous justification of this phenomena can be found in \cite{BurgosII}.

\section{Equations of Motion}

\begin{figure}[!hbp]
\centering
\includegraphics[width=3.5in]{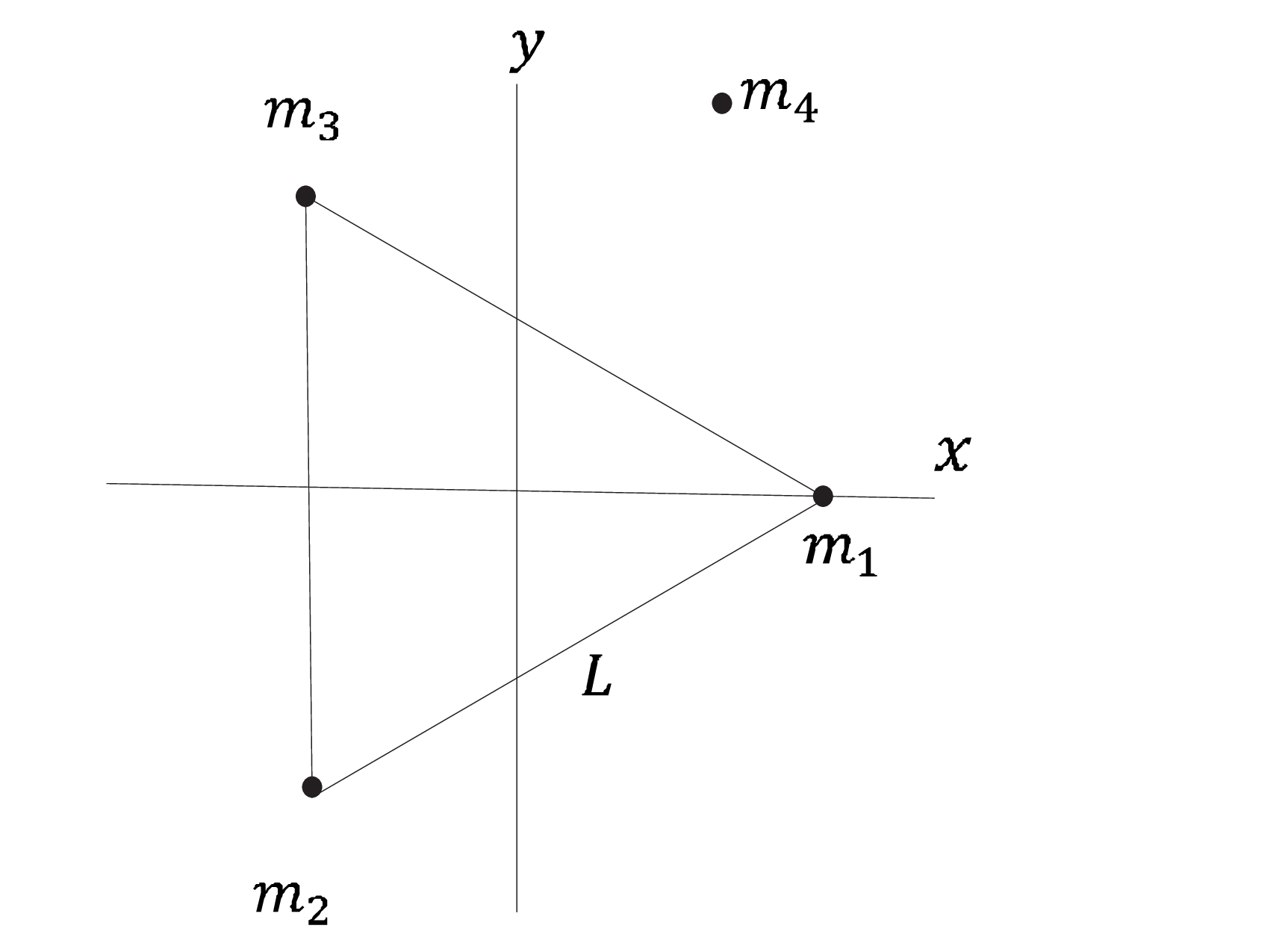}
\caption{The restricted four-body problem in a synodic system\label{triangle}}
\end{figure}

Consider three point masses, called $\textit{primaries}$, moving in
circular periodic orbits around their center of mass under their
mutual Newtonian gravitational attraction,  forming an equilateral
triangle configuration. A third massless particle moving in the same
plane is acted upon the attraction of the primaries. The equations
of motion of the massless particle referred to a synodic frame with
the same origin, where the primaries remain fixed, are:
\begin{eqnarray}
\bar{x}''-2n\bar{y}'-n^2\bar{x}&=&-k^2\sum_{i=1}^{3}m_{i}\frac{(\bar{x}-\bar{x_{i}})}{\rho_{i}^3}\nonumber\\
\bar{y}''+2n\bar{x}'-n^2\bar{y}&=&-k^2\sum_{i=1}^{3}m_{i}\frac{(\bar{y}-\bar{y_{i}})}{\rho_{i}^3} \label{sistemaconunidades}
\end{eqnarray}
where $k^2$ is the gravitational constant, $n$ is the mean motion,
$\rho_{i}^{2}=(\bar{x}-\bar{x}_{i})^2+(\bar{y}-\bar{y}_{i})^2$ is the distance of the massless particle to the primaries,
$\bar{x}_{i}$, $\bar{y}_{i}$
are the vertices of equilateral triangle  formed by the primaries,
and ($'$)  denotes derivative with respect to time $t^{*}$. We
choose the orientation of the triangle of masses such that $m_1$
lies along the positive $x$--axis and $m_2$, $m_3$ are located
symmetrically with respect to the same axis, see
figure~\ref{triangle}.

The equations of motion can be recast in dimensionless form as follows:
Let $L$ denote the length of triangle formed by the primaries,  $x=\bar{x}/L$,
$y=\bar{y}/L$, $x_i=\bar{x}_i/L$, $y_i=\bar{y}_i/L$, for $i=1,2,3$;
$M=m_{1}+m_{2}+m_{3}$ the total mass, and  $t=nt^{*}$. Then the equations (\ref{sistemaconunidades}) become
\begin{eqnarray}
\ddot{x}-2\dot{y}-x &=&-\sum_{i=1}^{3}\mu_{i}\frac{(x-x_{i})}{r_{i}^3}\nonumber\\
\ddot{y}+2\dot{x}-y &=&-\sum_{i=1}^{3}\mu_{i}\frac{(y-y_{i})}{r_{i}^3} \label{sistemasinunidades}
\end{eqnarray}
where we have used Kepler's third law: $k^2M=n^2L^3$, and the dot
($\dot{}$) represents derivatives with respect to the dimensionless
time $t$ and $r_{i}^2=(x-x_{i})^2+(y-y_{i})^2$. The system
(\ref{sistemasinunidades}) will be defined if we know the vertices
of triangle for each value of the masses. In this paper we suppose
$\mu:=\mu_{3}=\mu_{2}$ then $\mu_{1}=1-2\mu$, it's not hard to prove
that the vertices of triangle are given as function of the mass
parameter $\mu$ by $x_{1}=\sqrt{3}\mu$, $y_{1}=0$,
$x_{2}=-\frac{\sqrt{3}(1-2\mu)}{2}$, $y_{2}=-\frac{1}{2}$,
$x_{3}=-\frac{\sqrt{3}(1-2\mu)}{2}$, $y_{3}=\frac{1}{2}$. The system
(\ref{sistemasinunidades}) can be written succinctly as

\begin{eqnarray}
\ddot{x}-2\dot{y}&=&\Omega_{x} \label{sistemastandar}\\
\ddot{y}+2\dot{x}&=&\Omega_{y}
\end{eqnarray}
where
\begin{displaymath}
\label{omega}\Omega(x,y,\mu):=\frac{1}{2}(x^2+y^2)+\sum_{i=1}^{3}\frac{\mu_{i}}{r_{i}}.
\end{displaymath}
is the effective potential function.

 In the Restricted four--body problem (R4BP) there are three limiting cases:
\begin{enumerate}
\item If $\mu=0$, we obtain the rotating Kepler's problem, with $m_{1}=1$ at the origin of coordinates.
\item If $\mu=1/2$, we obtain the circular restricted three body problem, with two equal masses $m_{2}=m_{3}=1/2$.
\item If $\mu=1/3$, we obtain the symmetric case with three masses equal to $1/3$.
\end{enumerate}

It will be useful to write the system (\ref{sistemastandar}) using
complex notation. Let $z=x+\textit{i}y$, then
\begin{equation}
\label{sistemacomplejo}\ddot{z}+2\textit{i}\dot{z}=2\frac{\partial\Omega}{\partial\bar{z}}
\end{equation}
with
\begin{displaymath}
\Omega(z,\bar{z},\mu)=\frac{1}{2}\vert
z\vert^2+U(z,\bar{z},\mu)
\end{displaymath}
where the gravitational potential is
\begin{displaymath}
U(z,\bar{z},\mu)=\sum_{i=1}^{3}\frac{\mu_{i}}{\vert z-z_{i}\vert}
\end{displaymath}
and $r_{i}=\vert z-z_{i}\vert$, $i=1,2,3.$ are the distances to the
primaries. System (\ref{sistemacomplejo}) has the Jacobian first
integral
\begin{displaymath}
2\Omega(z,\bar{z},\mu)-\vert\dot{z}\vert^{2}=C
\end{displaymath}

If we define $P=p_{x}+\textit{i}p_{y}$, the conjugate momenta of
$z$, then system (\ref{sistemastandar}) can be recast as  a
Hamiltonian system  with Hamiltonian
\begin{eqnarray}
H&=&\frac{1}{2}\vert P\vert^2+Im(z\overline{P})-U(z,\bar{z},\mu)\nonumber\\
&=&\frac{1}{2}(p^{2}_{x}+p^{2}_{y})+(yp_{x}-xp_{y})-U(x,y,\mu).\label{hamiltoniano}
\end{eqnarray}
The relationship with the Jacobian integral is $H=-C/2$.
The phase space of (\ref{hamiltoniano}) is defined as
\begin{displaymath}
\Delta=\{(z,P)\in\mathbb{C}\times\mathbb{C}\vert z\ne z_{i},
i=1,2,3\},
\end{displaymath}
with collisions occurring at $z=z_{i}$, $i=1,2,3$. In the restricted three-body problem there exist five equilibrium
points for all values of the masses of the primaries but in this restricted
four-body problem the number of equilibrium
points depends on the particular values of the masses. A complete discussion of the equilibrium points and
bifurcations can be found in \cite{Del}, \cite{MeyerCC}, \cite{Lea},  \cite{PapaII}, \cite{Simo}.

\section{Regularization}
Where the solutions of the R4BP have binary collisions with any of the primaries the Hamiltonian (\ref{hamiltoniano}) is not defined for these solutions, so we have to remove such singularities in the system. The so called \textit{regularization process} is a technique that enable us to remove singularities of differential equations, therefore we want to apply this technique to the R4BP to study the system when the solutions are near to collision with the primaries. The \textit{regularization process} is a standar procedure and it can be found in \cite{Sz} and \cite{Gia}, however, we are going to explain it briefly in the present problem.\\ \\
First, we perform a translation from the center of mass, namely  $z=u+\sqrt{3}\mu-\sqrt{3}/2$, where $u=x_{2}+\textit{i}x_{2}$. The positions of the primaries in these new coordinates become $u_{1}=\frac{\sqrt{3}}{2}$, $u_{2}=-\frac{\textit{i}}{2}$, $u_{3}=\frac{\textit{i}}{2}$. In these coordinates the Hamiltonian is written as
\begin{equation}
\label{uhamiltonian}H=\frac{1}{2}\vert
U\vert^2+Im((u+\sqrt{3}\mu-\sqrt{3}/2)\overline{U})-V(u,\bar{u},\mu)
\end{equation}
where $U=P$ is the complex conjugate momenta de $u$.
We will denote by $f^{*}(w)$ the derivative with respect to complex variable $w$, $f(w)$ represents a complex valued analytic map. The following lemma shows how to complete a point transformation given by an analytic function $u=f(w)$
to a canonical transformation. We will chose later the mapping $f(w)$ to eliminate the singularities due to the primaries.

\begin{lem} Let $u=f(w)$ be a transformation point, then the transformation of the conjugate momenta $U=W/\overline{f^{*}(w)}$ yields a canonical transformation whenever $f^{*}(w)\ne0$
\end{lem}
$\textit{Proof}$: The mapping $(u,U)\rightarrow(w,W)$ is canonical if
$$ Re(\overline{U}du)=Re(\overline{W}dw)
$$ But
$$\overline{U}du=\bar{U}f^{*}(u)du=\frac{\overline{W}}{f^{*}(u)}f^{*}(u)du=\overline{W}dw
$$ form which the result follows.
\begin{flushright}
$\Box$
\end{flushright}
The regularization process starts with a canonical transformation followed by a re--parametrization of time on a fixed energy level $H=-C/2$. Let $f(w)$ be as above, this transformation must satisfy the hypothesis of the previous lemma. We perform the following scaling of time
\begin{equation}
\label{escalamientotemporal} \frac{d\tau}{dt}=\frac{1}{\vert f^{*}(w)\vert^{2}}
\end{equation}
Now we need to transform the Hamiltonian (\ref{uhamiltonian}) to the new variables
\begin{equation}
\label{hamiltonianow}H=\frac{1}{2}\frac{\vert W\vert^{2}}{\vert
f^{*}(w)\vert^{2}}+Im\left((f(w)+\sqrt{3}\mu-\sqrt{3}/2)\frac{\overline{W}}{f^{*}(w)}\right)-V(f(w),\overline{f(w)},\mu)
\end{equation}
Next perform Poincare's trick to re-parametrize solutions according to (\ref{escalamientotemporal})
$\overline{H}=\vert f^{*}(w)\vert^{2}(H+C/2)$. Observe that the energy level $H=-C/2$ is carried on to the level $\overline{H}=0$, explicitly
\begin{equation}
\label{hamiltonianobarra}\overline{H}=\frac{1}{2}\vert
W\vert^{2}+Im((f(w)+\sqrt{3}\mu-\sqrt{3}/2)\overline{f^{*}(w)W})-\vert
f^{*}(w)\vert^{2}V(w,\bar{w},\mu)
\end{equation}
\begin{displaymath}
+\vert f^{*}(w)\vert^{2}(C/2)
\end{displaymath}
where $$V(w,\bar{w},\mu)=\sum_{i=1}^{3}\frac{\mu_{i}}{\vert
f(w)-w_{i}\vert}$$ and $w_{i}$ denotes the position of the primaries.
Now we must choose the transformation $f(w)$ according to the conditions mentioned above and keeping in mind the singularities that we want to remove. Note that if we want to remove a single collision with any of the primaries, we can apply the Levi--Civita transformation \cite{Sz} to remove such singularity, however, the interesting problem is to remove simultaneously two or more singularities. It is not hard to see that the equations of the R4BP have the property that if $z(t)$ is a solution (in complex notation) then $\bar{z}(-t)$ is also a solution, in other words, we have symmetry of the solutions with respect to the $x-$axis as in the R3BP. This symmetry of the equations tells us that a collision with the primary $m_{2}$ (respectively $m_{3}$) implies necessarily a collision with the primary $m_{3}$ (respectively $m_{2}$), therefore a simultaneous regularization with the primaries $m_{3}$ and $m_{3}$ is needed. If we want to perform a simultaneous regularization in this case, first, we must note the importance of making the regularized equations as simple as possible in order to simplify the calculations and to save time in the integration of the equations. Therefore, we choose a transformation $f(w)$ similar to the Birkhoff's transformation \cite{Birk}
\begin{equation}
\label{transformacion} u=f(w)=\frac{1}{2}\left(w-\frac{1}{4w}\right)
\end{equation}
It's easy to prove that $f(w)$ has the following properties
\begin{eqnarray}
u_{i}&=&f(u_{i}),\quad i=2,3.\label{fixed}\\
f^{*}(w)&=&\frac{1}{2}\frac{(w-u_{2})(w-u_{3})}{w^{2}}\label{fprima}
\end{eqnarray}
 In particular
\begin{equation}
\vert f^{*}(w)\vert^{2}=\frac{1}{4}\frac{(\vert w-u_{2}\vert^{2})(\vert w-u_{3}\vert^{2})}{\vert w\vert^{4}}\label{deriv}
\end{equation}
and
\begin{equation}
f^{*}(u_{i})=0,\quad i=2,3.
\end{equation}
Observe that the positions of the primaries $m_{2}$ and $m_{3}$ remain fixed under the transformation. In fact, the following properties are the key to remove the singularities $w=w_{i}$, $i=2,3.$
$$f(w)-u_{i}=\frac{1}{2}\left(\frac{(w-u_{i})^{2}}{w}\right)$$ and $$f(w)-u_{1}=\frac{1}{2}\left(\frac{(w-a_{1})(w-a_{2})}{w}\right)$$ where $a_{1}=1+\sqrt{3}/2$, $a_{2}=-1+\sqrt{3}/2$. We must check that the Hamiltonian (\ref{hamiltonianobarra}) is free of singularities due to collisions with the primaries $m_{2}$ and $m_{3}$, observe that these points are contained in the term $\vert
f^{*}(w)\vert^{2}V(w,\bar{w},\mu)$, a straightforward calculation using (\ref{deriv}) shows
$$ \vert f^{*}(w)\vert^{2}V(w,\bar{w},\mu)=$$ $$\frac{1}{2\vert
w\vert^{3}}\left(\frac{(1-2\mu)\vert w-u_{2}\vert^{2}\vert
w-u_{3}\vert^{2}}{\vert w-a_{1}\vert \vert w-a_{2}\vert}+\mu(\vert w-u_{2}\vert^{2}+\vert
w-u_{3}\vert^{2})\right)
$$
Observe that we have removed the singularities due to the primaries $m_{2}$ and $m_{3}$, however, we have introduced new singularities though, $w=a_{1}$, $w=a_{2}$ and $w=0$. We are going to study these new singularities.
The origin of the new system $w=0$ is mapped under (\ref{transformacion}) to infinity in $u-$space, so it corresponds to escapes and it is not of interest for us. Let us analyze the remaining singularities $w=a_{1}$ and $w=a_{2}$. First we want to describe the number of pre-images of a point $u$ under the transformation $f(w)$, we need to solve the equation given by (\ref{transformacion}), or to find the roots of the quadratic polynomial $p(w)=w^{2}-2uw-1/4$. Note that given $u$, we have two roots or pre-images counting multiplicities. We recall the next proposition:
\begin{lem}
Let $w_{0}$ be a root of the polynomial $p(w)$, if $p'(w_{0})=p(w_{0})=0$ but $p''(w_{0})\ne0$ then $w_{0}$ is a root with multiplicity 2. If $p''(w_{0})=p'(w_{0})=p(w_{0})=0$ but $p'''(w_{0})\ne0$ then $w_{0}$ is a root with multiplicity 3 etc.
\end{lem}
It's easy to see that $p'(w)=2w-2u$ and $p''(w)=2$ then $p'(w)=0$ $\Leftrightarrow$ $w=u$. Now if we evaluate this root in the polynomial $p(w)$ we see that $p(u)=-u^{2}-1/4$ and $p(u)=0$ $\Leftrightarrow$ $w=u_{i}$ $i=2,3.$ therefore $p(u_{i})=p'(u_{i})=0$ but $p''(u_{i})\ne0$. In conclusion we have proved the following
\begin{prop} Let $u\in\mathbb{C}$ be a complex number, if $u=u_{i}$, $i=2,3.$ then the number of pre-images is one, if $u\ne u_{i}$, $i=1,2,3.$ then the number of pre-images is two.
\end{prop}
\begin{flushright}
$\Box$
\end{flushright}
This shows that the number of pre-images of the positions of the primaries $m_{2}$ and $m_{3}$ is exactly one.  Actually for $u=u_{i}$, $i=2,3$, we have
$$
p(w)=(w-u_{i})^{2}
$$
therefore the pre-images of each $u_{i}$, $i=2,3$ are they self. The pre-images of the primary $u_{1}=$ are exactly $a_{1}$ and $a_{2}$, then the new singularities correspond to the singularity $u_{1}$ in the $u$-space however we are not interested in removing this singularity, see figure (\ref{hillregions}). Instead, we have performed a global regularization of the singularities due to collisions with the primaries $m_{2}$ and $m_{3}$.
The phase space where the Hamiltonian (\ref{hamiltonianobarra}) is regular is given by
$$
\overline{\Delta}=\{(w,W)\in\mathbb{C}\times\mathbb{C}\vert w\notin\{0,a_{1},a_{2}\}, i=1,2,3\}
$$
Since the Hamiltonian (\ref{hamiltonianobarra}) contains only quadratic modulus, its partial derivatives are continuous throughout the region $\overline{\Delta}$.

\section{Hill's Regions of the Regularized Hamiltonian}
\begin{figure}[!hbp]
\begin{center}
\includegraphics[height=18.0cm,width=4.8in]{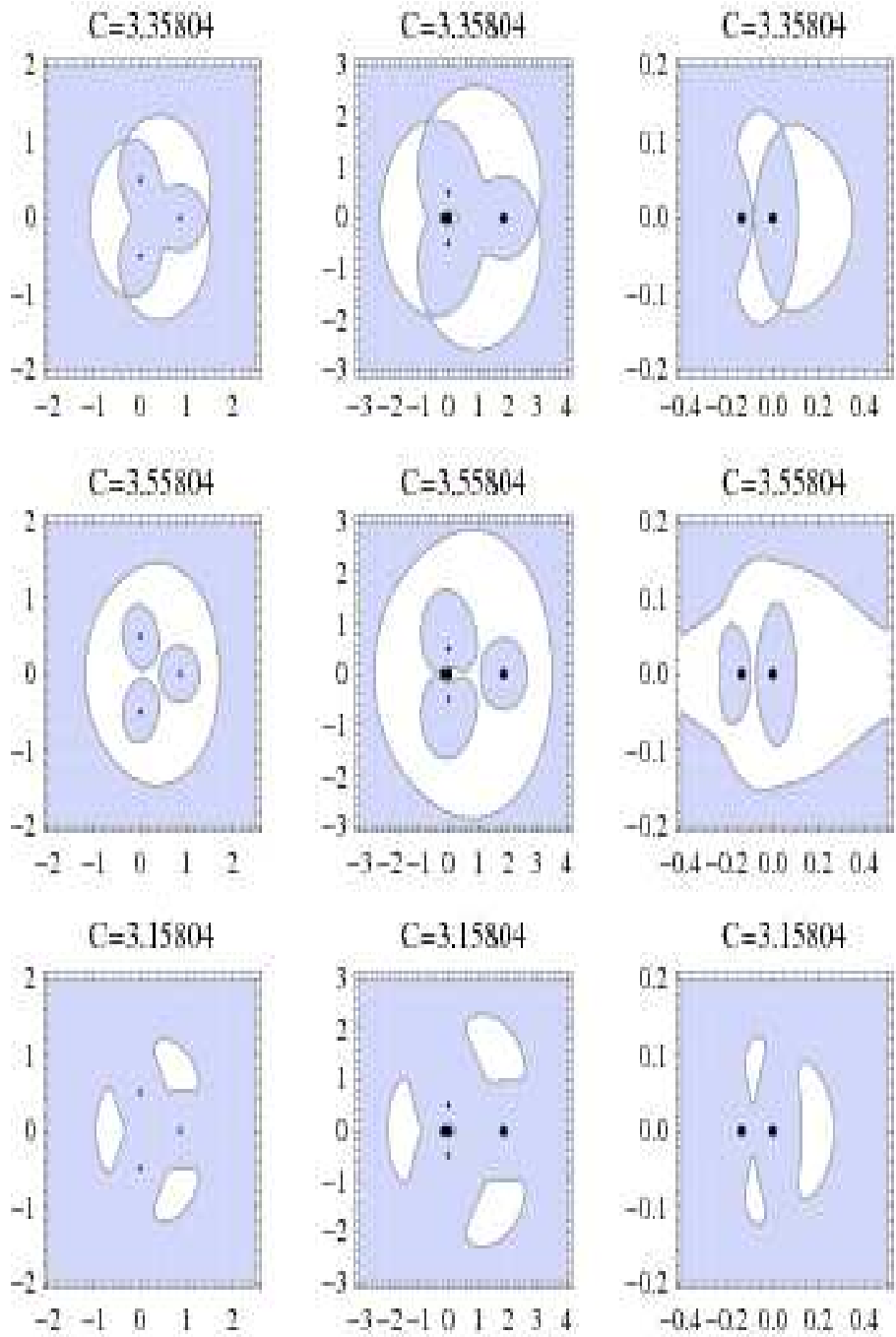}
\caption{Hill's regions (shaded areas) in the $u$-space (left column), in the $w$-space (center column) and magnifications of the origin in the $w$-space (right column). From top to bottom: $C=C_{1}$, $C=C_{1}+0.2$, $C=C_{1}-0.2$}\label{hillregions}
\end{center}
\end{figure}

The relation given by the first integral
$\vert\dot{z}\vert^{2}=2\Omega(z,\bar{z},\mu)-C$ implies
$2\Omega(z,\bar{z},\mu)-C\ge 0$ or $\Omega(z,\bar{z},\mu)\ge C/2$,
this inequality places a constraint on the position variable $z$ for
each values of $\mu$ and $C$, if $z$ satisfies this condition,
then there is a solution through that point $z$ for that values of
$C$ and $\mu$ (see for example \cite{MeyerHDS}). The sets where the inequality
$\Omega(z,\bar{z},\mu)\ge C/2$ holds are called the Hill's regions,
in regularized variables the former inequality becomes (see \cite{Sz})$$
\vert f^{*}(w)\vert^{2}(\Omega(w,\bar{w},\mu)-C/2)\ge0 $$ This
inequality defines the Hill's regions in regularized variables
whenever $w\ne0$, $w\ne a_{1}$ and $w\ne a_{2}$. Explicitly these regions
are defined by the expression \begin{equation}
\label{regioneshill}\frac{1}{2}\vert f^{*}(w)\vert^{2}\vert
f(w)+\sqrt{3}\mu-\frac{\sqrt{3}}{2}\vert^{2}+\vert
f^{*}(w)\vert^{2}V(w,\bar{w},\mu)-\vert
f^{*}(w)\vert^{2}\frac{C}{2}\ge0 \end{equation}
The next figure
shows the Hill's regions in the
$u-$space and the $w-$space for several values of the Jacobi
constant $C$ and for the equal masses case, i.e. $\mu=1/3$. At the reference value $C_{1}=3.35804$ there exist three critical points of the potential $\Omega$ in the $u$-space and the Hill's regions are very similar, however at the origin of the $w$-space new regions appear around the singularities $a_{1}$ and $0$, see figure (\ref{hillregions}). If we increase the value of $C_{1}$ the Hill's regions become disconnected in both spaces and the new regions around the singularities in the $w$-space increase their size, at this point it is clear the correspondences between the $u$ and $w$ spaces given by the transformation (\ref{transformacion}) discussed in the section (3). Finally, if we decrease the value of $C_{1}$ we find that the whole Hill's region is now connected, see figure (\ref{hillregions}). The positions of the primaries are marked by small circles and the singularities $w=0$, $w=a_{1}$ and $w=a_{2}$ are marked by black points.
\section{An application of the regularization process}
\begin{figure}[!hbp]
\begin{center}
\includegraphics[height=8.0cm, width=4.8in]{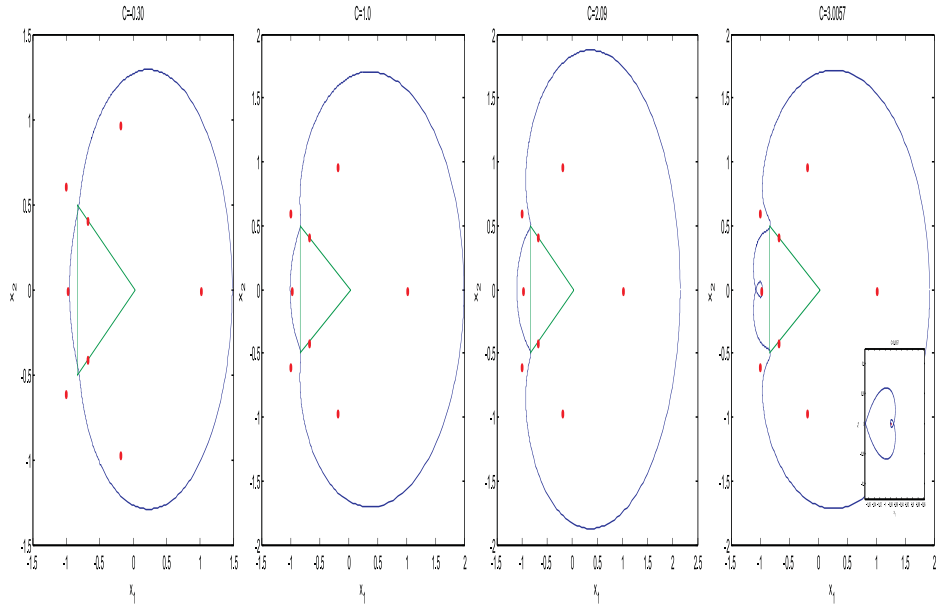}
\end{center}
\caption{The evolution of the second phase of family $f$ .\label{phasesf}}
\end{figure}

\begin{figure}[!hbp]
\begin{center}
\includegraphics[height=8.0cm, width=4.8in]{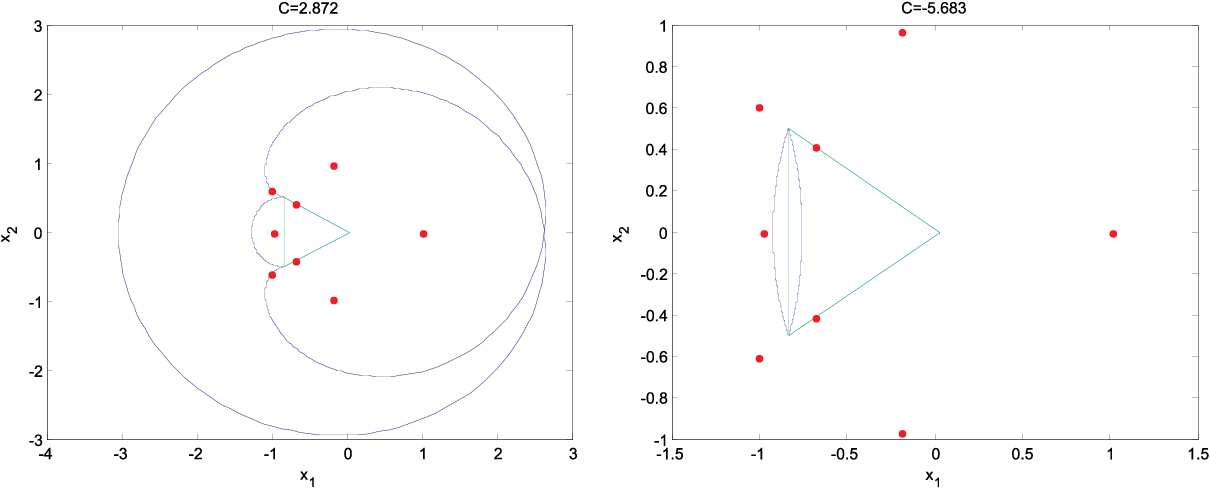}
\end{center}
\caption{Ejection--collision orbits of the family $j$, in the first phase (right) and in the second phase (left).\label{colisionesj}}
\end{figure}

\begin{figure}[!hbp]
\begin{center}
\includegraphics[height=8.0cm, width=4.8in]{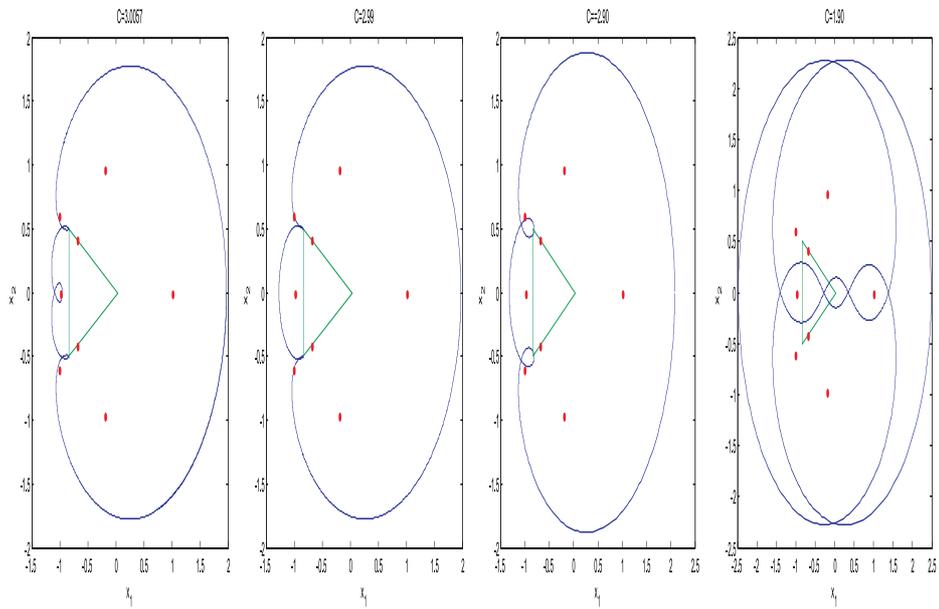}
\end{center}
\caption{The evolution of the first phase of family $r_{2}$ .\label{fasesr2}}
\end{figure}

\begin{figure}[!hbp]
\begin{center}
\includegraphics[height=8.0cm, width=4.8in]{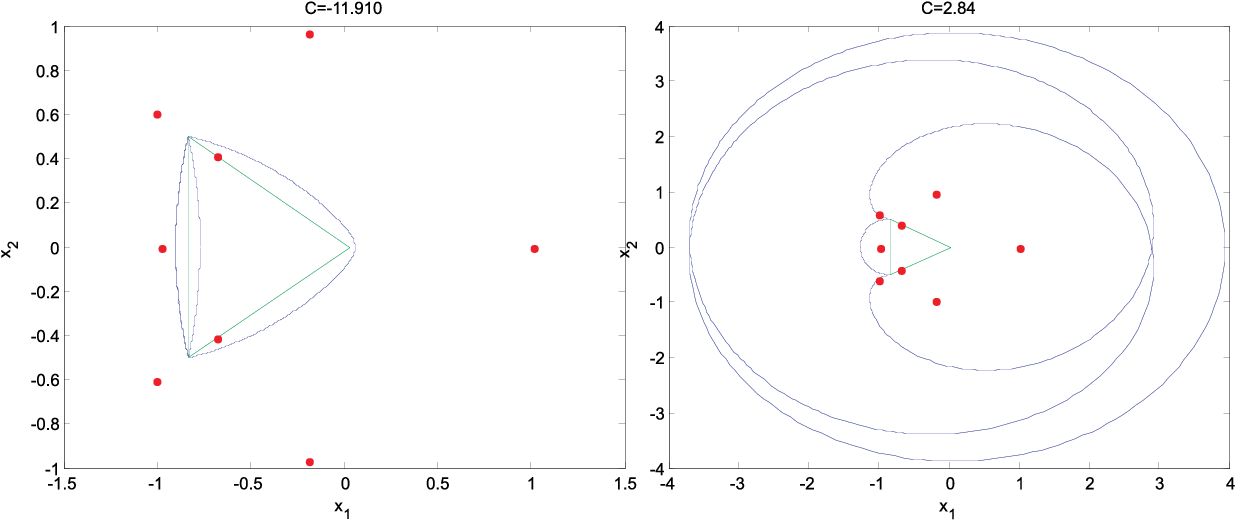}
\end{center}
\caption{Ejection--collision orbits of the family $j_{2}$, in the first phase (left) and in the second phase (right).\label{colisionesj2}}
\end{figure}

In this section we show an application of the regularized equations of the R4BP. We recall that Routh's criterion for linear stability of the Lagrangian configuration states that the masses of primaries must satisfy the inequality
\begin{displaymath}
\frac{m_{1}m_{2}+m_{2}m_{3}+m_{3}m_{1}}{m_{1}+m_{2}+m_{3}}<\frac{1}{27}.
\end{displaymath}
When the three masses are such that $m_{2}=m_{3}:=\mu$ and
$m_{1}+m_{2}+m_{3}=1$, the inequality is satisfied in the interval
$\mu\in[0,0.0190636...)$. In the paper \cite{Burgos} we can find a numerical exploration of families periodic orbits of the R4BP for values of the masses satisfying the Routh's criterion; in that work, there are nine families of periodic orbits and some of them contain ejection--collision orbits with the primaries $m_{1}$ and $m_{2}$, we are going to explain briefly how these orbits were obtained. Suppose that we have calculated a periodic orbit, this periodic orbit lies on a surface defined by the Jacobian first integral and therefore it has a well defined value of the constant $C$, if we use the analytical continuation method \cite{HenI}, \cite{Sz} we can follow the evolution of this orbits as the value of the constant $C$ varies continuously, in this evolution, the periodic orbit can reach collisions with any of the primaries, if we want to follow the orbit beyond these collisions we need to use regularized equations. When a ejection--collision orbit is reached, we say that the periodic orbit finishes one phase because after this collision the orbit changes its behavior, for instance from direct to retrograde. We refer to the reader to the references to see a complete discussion on families of periodic orbits. In the following we show some ejection--collision orbits in the R4BP and we explain where these orbits can be found.\\ \\
In the second phase of the family $f$, all of the orbits are near to collision with the primaries $m_{2}$ and $m_{3}$, however, these collisions are never reached but the regularized equations are needed in the numerical calculations to follow the evolution of this phase and to state the ``Blue sky catastrophe" termination of this family, see figure~ \ref{phasesf}.
In the family $j$ we find two ejection--collision orbits, one of them is at the beginning of the first phase and the second one is at the end of the second phase, see figure~ \ref{colisionesj}.
In the evolution of the first phase of the family $r_{2}$, we find similar orbits to the family $f$, but in this case a ejection--collision orbit appears before the ``Blue sky catastrophe"  termination of this family, see figure~ \ref{fasesr2}.
Finally, in the family $j_{2}$ we find two ejection--collision orbits, in the figure~ \ref{colisionesj2} we can see that these orbits are very similar to the ones of the family $j$.

\newpage
\bibliographystyle{plain}

\bibliography{biblio-u1}

Departamento de Matem\'aticas
UAM--Iztapalapa. Av. San Rafael Atlixco 186, Col. Vicentina, C.P.
09340, M\'exico, D.F. e--mail: jbg84@xanum.uam.mx,
\end{document}